\documentclass[12pt,twoside]{article}
\usepackage{amsmath,amssymb,array,latexsym}
\newtheorem{thm}{Theorem}
\newtheorem{lem}{Lemma}
\newtheorem{cor}{Corollary}

\def\ve{\varepsilon}
\def\notin{\epsilon \hspace{-0.37em}/}
\def\card{{\rm Card}}

\begin{document}

\author{Rafa{\l} Lata{\l}a and Joel Zinn}
\date{Warsaw University and Texas A\&M University}
\title{Necessary and Sufficient Conditions for the Strong Law of 
  Large Numbers for U-statistics
  }

\maketitle
\markboth{R.\ LATA{\L}A AND J.\ ZINN}{SLLN FOR U-STATISTICS}

\begin{abstract}
Under some mild regularity on the normalizing sequence, we obtain 
necessary and sufficient conditions for the Strong Law of Large 
Numbers for (symmetrized) U-statistics. We also obtain nasc's for 
the a.s. convergence of series of an analogous form. 
\end{abstract}

\section{Introduction.}

The general question addressed in this paper is that of necessary 
and sufficient conditions for 

\[
   \frac{1}{\gamma_{n}}\sum_{{\bf i}\in I_{n}}\ve_{{\bf i}} 
h(X_{{\bf i}})\to 0, 
\text{ a.s. },\]
where $I_{n}=\{{\bf i}=
(i_{i},i_{2},\ldots,i_{d}):1\leq i_{1}<i_{2}<\ldots<i_{d}\leq n\}$, 
$\{X_j\}_{j=1}^\infty$ is a sequence of iid r.v.'s,  
$X_{{\bf i}}=(X_{i_1},\cdots,X_{i_d})$. With no loss of generality
we may assume that $h$ is symmetric in its 
arguments. 

Further, as in \cite{CGZ} and in \cite{Zhang prod}, it is also important to 
consider the question of the almost sure convergence to zero of 
\[\frac{1}{\gamma_{n}}\max_{{\bf i}\in I_{n}}|h(X_{{\bf i}})|.\]
In fact, it is through the study of this problem that one 
is able to complete the characterization for the original question. 
 
Without the symmetrization by Rademachers, Hoeffding ([H]) in 1961  proved 
that for general $d$ and $\gamma_n=\binom{n}{d}$, mean zero is 
sufficient for the normalized sum above 
to go to zero almost surely. And, under a 
$p^{\hbox{th}}$ moment one has the a.s. 
convergence to zero with  $\gamma_n=n^{\frac{d}{p}}$ (\cite{S} when 
$0 < p< 1$, 
in the product case with mean zero \cite{T} for $1\le p< 2$ and in the case of 
general degenerate $h$ \cite{GZ} for $1<p<2$).

It is somewhat surprising that it took until the 90's to see that Hoeffding's 
sufficient condition was not necessary (\cite{GZ}). In the particular case 
in which $d=2$, $h(x,y)=xy$ and the variables are symmetric, 
necessary and sufficient conditons were 
given in (\cite{CGZ}) in 1995. This was later extended to $d\ge 3$ by Zhang 
(\cite{Zhang prod}). Very recently Zhang \cite{Zhang gen} 
obtained ``computable'' necessary and 
sufficient conditions in the case $d=2$ and, in general, found equivalent 
conditions in terms of a law of large numbers for modified maxima. Other 
related work is that of \cite{McConnell} in which the different indices go to 
infinity at their own pace and \cite{G} in which 
the variables in different coordinates can be based on different 
distributions.

In this paper we obtain nasc's for strong laws for `maxima' for general 
$d$. This likely would have enabled one to complete Zhang's program. 
However, we also found a more classical way of handling the reduction 
of the case of sums to the case of max's. 

The organization of the paper is as follows. In Section 2 we introduce 
the necessary notation and give the basic Lemmas. Now the form of our main 
Theorem is inductive. The reason we present the result in this form is that 
the conditions in the case $d>2$ are quite involved.
Because of the format of our Theorem we first present in Section 3, 
the case that the function, $h$, is 
the product of the coordinates. As mentioned earlier, this case received quite 
a bit attention,  culminating in Zhang's paper (\cite{Zhang prod}). In 
the first part of
Section 3 we show how the methods developed in this paper allow one to give a 
relatively simple, and perhaps transparent, proof of Zhang's result. 
We, then, prove the main result, namely, the nasc's for the Strong 
Law for symmetric U-statistics. Again, because of our inductive format, in 
order to clearly bring out the main idea's of our proof, 
we also give a simple proof of Zhang's result for 
the case $d=2$.

Finally in Section 4 we consider the question of 
convergence of multidimensional random series 
$\sum_{{\bf i}\in Z_{+}^{d}}h_{{\bf i}}(\tilde{X}_{{\bf i}})$. We obtain
necessary and sufficient conditions for a.s.\ convergence in the
case of nonnegative or symmetrized kernels. This generalizes the results
of \cite{KW1} (case $d=2$ and $h_{i,j}(x,y)=a_{i,j}xy$).

\section{Preliminaries and Basic Lemmas.}

Let us first introduce multiindex notation we will use in the paper:
\begin{itemize}
\item ${\bf i}=(i_{i},i_{2},\ldots,i_{d})$-multiindex of size $d$
\item $X_{{\bf i}}=(X_{i_{i}},X_{i_{2}},\ldots,X_{i_{d}})$, where $X_{j}$ is a
 sequence of i.i.d. random variables with values in some space $E$ and the 
common law $\mu$
\item $\tilde{X}_{{\bf i}}=(X^{(1)}_{i_{i}},X^{(2)}_{i_{2}},\ldots,
  X^{(d)}_{i_{d}})$,
 where $(X_{j}^{(k)})$, $k=1,\ldots,d$ are independent copies of $(X_{j})$,
\item $\ve_{{\bf i}}=\ve_{i_{1}}\ve_{i_{2}}\cdots\ve_{i_{d}}$, where $(\ve_{i})$ is
 a Rademacher sequence (i.e. a sequence of i.i.d.\ symmetric
 random variables taking on values $\pm 1$) independent of other random
  variables
\item $\tilde{\ve}_{{\bf
i}}=\ve_{i_{1}}^{(1)}\ve_{i_{2}}^{(2)}\cdots\ve_{i_{d}}^{(d)}$, 
where $(\ve_{i}^{(j)})$ is
 a doubly indexed Rademacher sequence independent of other random 
 variables
\item $\mu_{k}=\otimes_{i=1}^{k}\mu$ - product measure on $E^{k}$
\item for $I\subset\{1,2,\ldots,d\}$, by $E_{I}$ and $E_{I}^{'}$  we will 
    denote expectation with respect to $(X_{i}^{k})_{k\in I}$ and
    $(X_{i}^{k})_{k\notin I}$ respectively
\item ${\bf i}_{I}=(i_{k})_{k\in I}$ and $I^{'}=\{1,2,\ldots,d\}
      \setminus I$ for $I\subset\{1,2,\ldots,d\}$
\item $I_{n}=\{{\bf i}=(i_{i},i_{2},\ldots,i_{d}):1\leq i_{1}<i_{2}<\ldots<i_{d}\leq 
 n\}$,      
\item $C_{n}=\{{\bf i}=(i_{i},i_{2},\ldots,i_{d}):1\leq i_{1},i_{2},\ldots,i_{d}\leq 
 n\}$
\item $A^{I,x}=A^{x_{I}}=\{z\in E^I: \exists a\in A, a_I=x_I, a_{I^{\prime}}=z\}$ for $A\subset E^{d}, I\subset \{1,\ldots,d\}$.
\end{itemize}

\bigskip

  The results in this section were motivated by the difficulty 
  in computing quantities such as:
\[P(\max_{i,j\le n}h(X_i,Y_j)>t),\]
 where $\{X_i\}$ are independent random variables and $\{Y_i\}$ is 
 an independent copy, and $h$ is, say, symmetric in its arguments. 

 In the one-dimensional case, namely,  
 $P(\max_{i\leq n} \xi_i > t)$, where $\{\xi_i\}$ are independent 
 r.v.'s, we have the simple inequality
\begin{equation}
\label{d=1}
  \frac{1}{2}\min(\sum_i P(|\xi_i| > t),1)\leq P(\max_i|\xi_{i}|> t)
  \leq \min(\sum_i P(|\xi_{i}| > t),1).
\end{equation}

 If this type of inequality held for any dimension, the proofs and 
 results would look much the same as in dimension $1$. Here we give an 
 example to see the difference between the cases $d=1$ and $d>1$. 

  Consider the set in the unit square given by:
\[A=\{(x,y)\in [0,1]^2 : x<a,y<b \mbox{ or } x<b, y<a\}\]
  and assume that the $X_i, Y_{j}$ are iid uniformly distributed on $[0,1]$.
  By (\ref{d=1}) it easily follows that
  
 \[P(\max_{1\leq i,j \leq n}I_{A}(X_{i},Y_{j})>0)\sim \min(na,1)\min(nb,1),\]
 which is equivalent to $\sum_{i,j=1}^{n}P(I_{A}(X_{i},Y_{j})>0)\sim n^{2}ab$
 if and only if both $a$ and $b$ are of order $O(\frac{1}{n}).\Box$

\begin{lem}
\label{techdec}
  Suppose that the nonnegative functions $f_{{\bf i}}(x_{{\bf i}})$ satisfy the 
  following conditions
 \begin{equation}
 \label{ass1.1}
   f_{{\bf i}}(\tilde{X}_{{\bf i}})\leq 1  \mbox{\ a.s. for all {\bf i}}
 \end{equation}
 \begin{equation}
 \label{ass1.2}
   E_{I}\sum_{{\bf i}_{I}}f_{{\bf i}}(\tilde{X}_{{\bf i}})\leq 1 
   \mbox{\ a.s. for any $I\subset\{1,2,\ldots,d\}$, $0<\card(I)<d$ }
 \end{equation}
  Let $\tilde{m}_{1}=E\sum_{{\bf i}}f_{{\bf i}}(\tilde{X}_{{\bf i}})$, then
 \begin{equation}
 \label{in1.1}
   E(\sum_{{\bf i}}f_{{\bf i}}(\tilde{X}_{{\bf i}}))^{2}\leq \tilde{m}_{1}^{2}+
   (2^{d}-1)\tilde{m}_{1}
 \end{equation}  
  and
 \begin{equation}
 \label{in1.2}
   P(\sum_{{\bf i}}f_{{\bf i}}(\tilde{X}_{{\bf i}})\geq 
   \frac{1}{2}\tilde{m}_{1})\geq 
   2^{-d-2}\min(\tilde{m}_{1},1).
\end{equation}  
\end{lem}

\noindent{\bf Proof.} Let $S(d)$ denote the family of nonempty subsets of 
  $\{1,\ldots,d\}$
  and for a fixed $I\in S(d)$ and ${\bf i}$ let
\[\tilde{J}({\bf i},I)=\{{\bf j}: {\bf j}_{I}={\bf i}_{I} 
  \mbox{ and $j_{k}\neq i_{k}$ for all $k\notin I$}\}.\]
 Then we have by (\ref{ass1.1}) and (\ref{ass1.2})
\[E(\sum_{{\bf i}}f_{{\bf i}}(\tilde{X}_{{\bf i}}))^{2}\leq
  (E\sum_{{\bf i}}f_{{\bf i}}(X_{{\bf i}}))^{2}
  +\sum_{I\in S(d)}\sum_{{\bf i}}E_{I}E^{'}_{I}f_{{\bf i}}
  (\tilde{X}_{{\bf i}})E^{'}_{I}\sum_{j\in \tilde{J}({\bf i},I)}f_{{\bf j}}
  (\tilde{X}_{{\bf j}})\]
\[\leq m_{1}^{2}+\sum_{I\in S(d)}\sum_{{\bf i}}E_{I}E^{'}_{I}f_{{\bf i}}
    (\tilde{X}_{{\bf i}})= m_{1}^{2}+(2^{d}-1)m_{1}.\]    
The inequality (\ref{in1.2}) follows by (\ref{in1.1}) and the Paley-Zygmund
inequality.$\, \Box$
 
\medskip
 
The next Lemma is an undecoupled version of Lemma 1, the proof of it is similar
as of Lemma 1 and is omitted.
\medskip
\begin{lem}
\label{techund}
Suppose that the nonnegative functions $f_{{\bf i}}(x_{{\bf i}})$ satisfy the 
following conditions
\[
  f_{{\bf i}}(X_{{\bf i}})\leq 1  \mbox{\ a.s. for all {\bf i}}\]
and
\[E_{I}^{'}\sum_{j\in J({\bf i},I)}f_{{\bf j}}(X_{{\bf j}})\leq 1 \mbox{\ a.s. 
  for all {\bf i} and $I\subset\{1,2,\ldots,d\}$, $0<\card(I)<d$ },\]
where
\[J({\bf i},I)=\{{\bf j}:\{k:\exists_{l}\ i_{k}=j_{l}\}=I\}.\]  
Let $m_{1}=E\sum_{{\bf i}}f_{{\bf i}}(X_{{\bf i}})$, then
\begin{equation}
\label{in2.1}
  E(\sum_{{\bf i}}f_{{\bf i}}(X_{{\bf i}}))^{2}\leq m_{1}^{2}+
  (2^{d}-1)m_{1}
\end{equation}  
and
\begin{equation}
\label{in2.2}
  P(\sum_{{\bf i}}f_{{\bf i}}(X_{{\bf i}})\geq \frac{1}{2}m_{1})\geq 
  2^{-d-2}\min(m_{1},1).
\end{equation}  
\end{lem} 

In the rest of this paper we will refer to the next Corollary as the 
``Section Lemma''. 
 
\begin{cor}
 If the set $A\subset E^{d}$ satisfies the condition
 \[n^{d-l}\mu_{d-l}(A^{I,X_{I}})\leq 1 \mbox{\ a.s. for all $I\subset \{1,\ldots,d\}$
   with $0<\card(I)=l<d$}\]
 then 
 \[P(\exists_{{\bf i}\in C_{n}} \tilde{X}_{{\bf i}}\in A)\geq 2^{-d-2}
   \min(n^{d}\mu_{d}(A),1)\]
  and for $n\geq d$
  \[P(\exists_{{\bf i}\in I_{n}} X_{{\bf i}}\in A)\geq 2^{-d-2}d^{-d}
   \min(n^{d}\mu_{d}(A),1).\] 
\end{cor} 

\noindent{\bf Proof.} The first inequality follows immediately by Lemma \ref{techdec} applied
   to $f_{{\bf i}}=I_{A}$. To prove the second inequality we use Lemma 
   \ref{techund} and notice that
  \[\min(\binom{n}{d} \mu_{d}(A),1)\geq d^{-d}\min(n^{d}\mu_{d}(A),1).\, 
     \Box\]

\bigskip

\section{Strong Laws of Large Numbers}

 We will assume in this section that the sequence $\gamma_{n}$ satisfy the 
 following regularity conditions 
 \begin{equation}
 \label{assd1}
   \gamma_{n} \mbox{ is nondecreasing} 
 \end{equation}
 \begin{equation}
 \label{assd2}
   \gamma_{2n}\leq C\gamma_{n} \mbox{ for any $n$}
 \end{equation}
 \begin{equation}
 \label{assd3}
   \sum_{k\geq l}\frac{2^{dk}}{\gamma^{2}_{2^{k}}}\leq 
   C\frac{2^{dl}}{\gamma^{2}_{2^{l}}} \mbox{ for any $l=1,2,\ldots$}
 \end{equation}

As mentioned in the Introduction we first give a proof of Zhang's 
result 
\cite{Zhang prod} for the product case i.e. $h(x)=\prod_{i=1}^{d}x_{i}$ for
$x\in R^{d}$.  To state the SLLN in this case we need to define numbers $c_{n}$
by the formula
\[c_{n}=\min\{c>0: nE(\frac{X^{2}}{c^{2}}\wedge 1)\leq 1\}.\]

\begin{thm}
  Assume that $h(x)=\prod_{i=1}^dx_i$, and that the r.v.'s $X_{i}$ are 
  symmetric. Then, under the regularity assumptions
  (\ref{assd1})-(\ref{assd3}), the following are equivalent:

 \begin{equation}
 \label{prodcon}
   \frac{1}{\gamma_{n}}\sum_{{\bf i}\in I_{n}}h(X_{\bf i})= 
   \frac{1}{\gamma_{n}}\sum_{{\bf i}\in I_{n}}\prod_{r=1}^{d}X_{i_{r}}
   \rightarrow 0 \mbox{\ a.s.}
\end{equation}
\begin{equation}
\label{Zprod}
    \sum_{k=1}^{\infty} 2^{kl}
    P(\prod_{r=1}^{l}X_r^2>\frac{\gamma_{2^k}^2}{c_{2^{k}}^{2(d-l)}}, 
    \min_{r\le l}X_r^2>c_{2^k}^2)<\infty 
\text{ for all } 1\le l\le d. 
\end{equation}
\end{thm}
\noindent{\bf Proof.} We give only the proof of the necessity of the conditions
 (\ref{Zprod}). The 
 sufficiency can be proved as in the  Theorem \ref{Sum2Max}. Let 
  \[T^{(r)}_n=\sum_{i_r=1}^n{X_{i_r}^{(r)}}^2\]
  and
  \[ T^{(r)}_n(c)=\sum_{i_r=1}^n{X_{i_r}^{(r)}}^2\wedge c^2.\]

{\bf Step 1.} We first reduce to the sum of squares, i.e.\ we will show that 
 condition (\ref{prodcon}) implies 
 \begin{equation}
 \label{squares}
   \gamma_{n}^{-2}\sum_{{\bf i}\in I_{n}}\prod_{r=1}^{d} X_{i_{r}}^{2}
   \rightarrow 0 \mbox{ a.s.}
 \end{equation}
 By the symmetry of $X$ we have that $\gamma_{n}^{-1}\sum_{{\bf i}\in I_{n}}
 \prod_{r=1}^{d}\ve_{i_{r}}X_{i_{r}}\rightarrow 0$ a.s. Thus
 for a.a.\ sequences $(X_{i})$, the Walsh sums (i.e. the linear combinations
 of products of d Rademachers) 
 converge to 0 a.s. Hence, they converge in probability. This implies (by a 
 result of Bonami about hypercontractivity of Walshes \cite{B}) that for a.a. 
 sequences $(X_{i})$, 
 $\gamma_{n}^{-2}\sum_{{\bf i}\in I_{n}}\prod_{r=1}^{d} X_{i_{r}}^{2}
 \rightarrow 0$ and (\ref{squares}) is proved.

{\bf Step 2.} We now go to a diadic subsequence and then decouple.
  By the Borel-Cantelli Lemma, the condition (\ref{squares}) implies that
  \[\forall_{\ve>0}\ \sum_{k=1}^{\infty}
    P(\sum_{{\bf i}\in I_{2^{k-1}}}\prod_{r=1}^{d} X_{i_{r}}^{2}\geq
    \ve\gamma_{2^{k}}^{2})<\infty. \]
 Now let us notice that 
 $I_{2^k}\supseteq \{{\bf i}\in I_{2^k}: (r-1)2^{k-l}< i_r\le r2^{k-l}\}$ 
 if $l$ is such that $2^l\ge d$. Moreover the random variables in these blocks 
 are independent of the other blocks, thus we obtain
 \[\forall_{\ve>0}\  \sum_{k=l+1}^{\infty}
    P(\sum_{{\bf i}\in C_{2^{k-l-1}}}\prod_{r=1}^{d} (X_{i_{r}}^{(r)})^{2}\geq
    \ve\gamma_{2^{k}}^{2})<\infty. \]
 Hence, using the regularity assumption (\ref{assd2})
 \begin{equation}
 \label{sumdec}
 \forall_{\ve>0}\
   \sum_{k=1}^{\infty}P(\prod_{r=1}^{d}T_{2^{k}}^{(r)}\geq\ve\gamma_{2^{k}}^{2})
   <\infty.
 \end{equation}

{\bf Step 3.} At this point we use 1-dimensional case of Lemma 1. We apply 
  it to
  \[c_{n}^{-2}T_{n}^{(r)}(c_{n})=\sum_{j=1}^{n}
    \frac{(X_{j}^{(r)})^{2}}{c_{n}^{2}}\wedge 1\]
  and notice that $Ec_{n}^{-2}T_{n}^{(r)}(c_{n})=1$ by the definition of $c_{n}$.
  We get that 
 \[P(T_{n}^{(r)}(c_{n})\geq \frac{1}{2}c_{n}^{2})
   =P(c_{n}^{-2}T_{n}^{(r)}(c_{n})\geq \frac{1}{2}Ec_{n}^{-2}T_{n}^{(r)}(c_{n}))
   \geq \frac{1}{8}.\]  
  Hence
  \[P(\prod_{r=l+1}^{d}T_{n}^{(r)}\geq \frac{c_{n}^{2(d-l)}}{2^{d-l}})\geq
    P(\prod_{r=l+1}^{d}T_{n}^{(r)}(c_{n})\geq \frac{c_{n}^{2(d-l)}}{2^{d-l}})
    \geq (\frac{1}{8})^{d-l}\]
  and
  \[P(\prod_{r=1}^{d} T_{n}^{(r)}\geq 2^{l-d}\gamma_{2^{k}}^{2})\geq
    (\frac{1}{8})^{d-l}P(\prod_{r=1}^{l} T_{n}^{(r)}\geq 
     \frac{\gamma_{2^{k}}^{2}}{c_{2^{k}}^{2(d-l)}}).\]  
  Thus condition (\ref{sumdec}) yields
  \begin{equation}
  \label{maxprod}
   \sum_{k=1}^{\infty} P(\max_{i_1,\ldots,i_l\leq 2^{k}}
   \prod_{r=1}^{l} ({X_{i_r}^{(r)}})^2
   > \frac{\gamma^2_{2_k}}{c_{2^k}^{2(d-l)}}) < \infty
  \end{equation}

\noindent Now, here is the {\bf main point}. 

{\bf Step 4.} At this point we need to replace the max inside the 
probability with $2^{kl}$ outside the probability. To do this 
we use the Section Lemma (Corollary 1). 

To get small sections there are a 
variety of choices. To obtain Zhang's result, we reduce the 
probabilities even further by intersecting the sets in the following manner. 
\[\sum_{k=1}^{\infty} P(\max_{i_1,\ldots,i_l\le 2^k}\prod_{r=1}^{l} 
 ({X_{i_r}^{(r)}})^2 I_{\{(X_{i_r}^{(r)})^2>c^2_{2^k}\}}
  > \frac{\gamma^2_{2_k}}{c_{2^k}^{2(d-l)}}) < \infty \]
To see why we have small sections, just note that 
\[P(X^2 > c^2_{2^k})\le \frac{E(X^2\wedge c^2_{2^k})}{c^2_{2^k}}=
  \frac{1}{2^k}.\]
Now we just use the Section lemma to get
\[\sum_{k=1}^{\infty} 2^{kl} P(\prod_{r=1}^l X_{r}^2 I_{\{X_{r}^2>c^2_{2^k}\}}
  > \frac{\gamma^2_{2_k}}{c_{2^k}^{2(d-l)}}) < \infty\]
Or, equivalently,
\[\sum_{k=1}^{\infty} 2^{kl}P(\prod_{r=1}^l 
  X_{r}^2 > \frac{\gamma^2_{2_k}}{ c_{2^k}^{2(d-l)}}, 
  \min_{1\le r\le l}X_{r}^2>c^2_{2^k}) < \infty,\] 
which yields (\ref{Zprod}). $\Box$

\medskip

 In Theorem \ref{Sum2Max} 
 we reduce the SLLN for symmetric or nonnegative kernels to a SLLN for  
 ``modified maxima''. To see 
  what this means consider the case $d=2$. Then,
  \begin{eqnarray*}
    A_{k,2}=\{(x,y)\in E^2:& h^2(x,y)\le \gamma_{2^k}^2,
    2^kE_{Y}h^2I_{h^2\le \gamma_{2^k}}(x,Y)\le  \gamma_{2^k}^2,\\
    &2^kE_{X}h^2I_{h^2\le \gamma_{2^k}}(X,y)\le  \gamma_{2^k}^2\}.
  \end{eqnarray*}  
 So that 
 \[\{\exists {\bf i}\in C_{2^k}, {\tilde X}_{\bf i}\notin A_{k,2}\}=
   \{\max_{{\bf i}\in C_{2^k}}\varphi(\tilde{X}_{\bf i})>\gamma_{2^{k}}^{2}\},\]
 where
 \[\varphi(x,y)=h^{2}(x,y)\vee 2^kE_{Y}h^2I_{h^2\le \gamma_{2^k}}(x,Y)
   \vee 2^kE_{X}h^2I_{h^2\le \gamma_{2^k}}(X,y).\]  

In \cite{Zhang gen} Zhang, using different methods, also reduced the probem to 
``modified maxima''. 
We continue in Theorem \ref{Max} to find nasc's for the SLLN for the maximum, 
which, hence, could also be used to complete Zhang's program. \par

\medskip
 For a measurable
 function $h$ on $E^{d}$ which is symmetric with respect to  permutations 
of the variables,  we define for $k=1,2,\ldots$
\[A_{k,1}=\{x\in E^{d}: h^{2}(x)\leq \gamma_{2^{k}}^{2}\}\]
 and for $l=1,\ldots,d-1$
\begin{eqnarray*}\lefteqn{A_{k,l+1}}\\
& &=\{x\in A_{k,l}: 2^{kl}E_{I}h^{2}I_{A_{k,l}}(x)\leq
   \gamma_{2^{k}}^{2} \mbox{ for all $I\subset\{1,2,\ldots d\}$, $\card(I)=l$}\}
   .\end{eqnarray*}

\begin{thm}\label{Sum2Max}
Suppose that assumptions (\ref{assd1})-(\ref{assd3}) are satisfied and the sets
$A_{k,l}$ are defined as above. Then the 
following conditions are equivalent:

\begin{equation}
\label{A}
   \frac{1}{\gamma_{n}}\sum_{{\bf i}\in I_{n}}\ve_{{\bf i}}h(X_{{\bf i}})
  \rightarrow 0 \mbox{\ a.s.}
\end{equation}   
\begin{equation}
\label{Apr}
   \frac{1}{\gamma_{n}}\sum_{{\bf i}\in C_{n}}\tilde{\ve}_{{\bf i}}h(\tilde{X}_{{\bf i}})
  \rightarrow 0 \mbox{\ a.s.}
\end{equation} 
\begin{equation}
\label{B} 
\frac{1}{\gamma^{2}_{n}}\sum_{{\bf i}\in I_{n}}h^{2}(X_{{\bf i}})
  \rightarrow 0 \mbox{\ a.s.} 
\end{equation}  
\begin{equation}
\label{Bpr} 
\frac{1}{\gamma^{2}_{n}}\sum_{{\bf i}\in C_{n}}h^{2}(\tilde{X}_{{\bf i}})
  \rightarrow 0 \mbox{\ a.s.} 
\end{equation}  
\begin{equation}
\label{C} 
\sum_{k=1}^{\infty}P(\exists_{{\bf i}\in I_{2^{k}}}\ X_{{\bf i}}\notin A_{k,d})
  <\infty
\end{equation}
\begin{equation}
\label{Cpr} 
\sum_{k=1}^{\infty}P(\exists_{{\bf i}\in C_{2^{k}}}\ 
  \tilde{X}_{{\bf i}}\notin A_{k,d})<\infty
\end{equation}  
\end{thm}

\noindent{\bf Proof.} (\ref{A})$\Rightarrow$(\ref{B}) and 
   (\ref{Apr})$\Rightarrow$(\ref{Bpr}) - proofs of these implications are the 
   same as in Proposition 4.7 in \cite{CGZ} (see also Step 1 in the
   proof of Theorem 1)
   
(\ref{B})$\Rightarrow$(\ref{Bpr}) Let $l$ be such that $2^{l}\geq d$. By the 
 regularity of $\gamma_{n}$ (\ref{assd1}),(\ref{assd2}) and the 
 Borel-Cantelli Lemma, (\ref{B}) and (\ref{Bpr}) are equivalent, respectively,
 to
 \begin{equation}
 \label{Beq}
  \sum_{k=1}^{\infty}P(\sum_{{\bf i}\in I_{2^{k}}}h^{2}(X_{{\bf i}})\geq 
  \varepsilon \gamma_{2^{k}}^{2})< \infty \mbox{ for all $\varepsilon>0$}
 \end{equation}
 and
 \begin{equation}
 \label{Bpreq}
  \sum_{k=l+1}^{\infty}P(\sum_{{\bf i}\in C_{2^{k-l}}}
   h^{2}(\tilde{X}_{{\bf i}})\geq 
  \varepsilon \gamma_{2^{k}}^{2})< \infty \mbox{ for all $\varepsilon>0$}.
 \end{equation} 
 Let
 \[D_{k}=\{{\bf i}: (m-1)2^{k-l} < i_{m}\leq m2^{k-l} 
   \mbox{ for $m=1,\ldots,d$}\},\]
 then  for $k\geq l$ we get
 \[P(\sum_{{\bf i}\in I_{2^{k}}}h^{2}(X_{{\bf i}})\geq 
  \varepsilon \gamma_{2^{k}}^{2})\geq P(\sum_{{\bf i}\in D_{k}}
  h^{2}(X_{{\bf i}})\geq \varepsilon \gamma_{2^{k}}^{2})\] 
 \[=P(\sum_{{\bf i}\in C_{2^{k-l}}}h^{2}(\tilde{X}_{{\bf i}})\geq 
  \varepsilon \gamma_{2^{k}}^{2})\] 
  and (\ref{Beq}) implies (\ref{Bpreq}).
  
(\ref{B})$\Rightarrow$(\ref{C}) We will prove by induction that for 
 $l\leq d$
 \begin{equation}
 \label{ind}
  \sum_{k=1}^{\infty}P(\exists_{{\bf i}\in I_{2^{k}}}\ X_{{\bf i}}
   \notin A_{k,l})<\infty
 \end{equation} 
 For $l=1$ (\ref{ind}) is $\sum_{k=1}^{\infty}P(\exists_{{\bf i}\in I_{2^{k}}}\
 h^{2}(X_{{\bf i}})> \gamma_{2^{k}}^{2})<\infty$ and follows easily by the 
 Borel-Cantelli lemma.
 Assume that (\ref{ind}) holds for $l\leq d-1$. To show it for $l+1$ it is
 enough to prove that for any $I$ with $\card(I)=l$
 \begin{equation}
 \label{ind2}
  \sum_{k=1}^{\infty}P(\exists_{{\bf i}_{I^{'}}\in I_{2^{k}}}\ 2^{kl}E_{I}
  h^{2}I_{A_{k,l}}(X_{{\bf i}})> \gamma_{2^{k}}^{2})<\infty
 \end{equation}
 By the symmetry of the kernel $h$ we may and will assume that
 $I=\{1,\ldots,l\}$.  From (\ref{B}) it follows that 
 \[\frac{1}{\gamma_{2^{k}}^{2}}\sum_{{\bf i}\in I_{2^{k}}}
   h^{2}I_{A_{k,l}}(X_{{\bf i}})\rightarrow 0 \mbox{ a.s.}.\]
 By the regularity of $\gamma_{2^{k}}$ (\ref{assd2}) and the Borel-Cantelli
 Lemma we get that
 \[\sum_{k=1}^{\infty}P(\sum_{{\bf i}\in I_{2^{k+1}}}h^{2}I_{A_{k,l}}
   (X_{{\bf i}})\geq \frac{1}{2}\gamma_{2^{k}}^{2})<\infty.\]
 But
 \[P_{I}(\sum_{{\bf i}\in I_{2^{k+1}}}h^{2}I_{A_{k,l}}
   (X_{{\bf i}})\geq \frac{1}{2}\gamma_{2^{k}}^{2})\geq
  P_{I}(\max_{{\bf i}_{I^{'}}\in J_{2^{k}}} \sum_{{\bf i}_{I}\in I_{2^{k}}}
   h^{2}I_{A_{k,l}}(X_{{\bf i}})\geq \frac{1}{2}\gamma_{2^{k}}^{2})\]
 \[\geq \max_{{\bf i}_{I^{'}}\in J_{2^{k}}}P_{I}( \sum_{{\bf i}_{I}\in I_{2^{k}}}
   h^{2}I_{A_{k,l}}(X_{{\bf i}})\geq \frac{1}{2}\gamma_{2^{k}}^{2}),\]
 where 
 \[J_{2^{k}}=\{(i_{1},\ldots,i_{d-l}): 2^{k}<i_{1}<i_{2}<\ldots<i_{d-l}\leq
   2^{k+1}\}.\]  
 Let us notice that by the definition of $A_{k,l}$ we have for any $J\subset I$
 with $\card(J)=m<l$
\[2^{km}E_{J} h^{2}I_{A_{k,l}}(X_{{\bf i}})\leq  \gamma_{2^{k}}^{2}.\]
Therefore by Lemma \ref{techund} we get that 
\[\max_{{\bf i}_{I^{'}}\in J_{2^{k}}}P_{I}( \sum_{{\bf i}_{I}\in I_{2^{k}}}
   h^{2}I_{A_{k,l}}(X_{{\bf i}})\geq \frac{1}{2}\gamma_{2^{k}}^{2})
   \geq 2^{-l-2},\]
 if $\max_{{\bf i}_{I^{'}}\in J_{2^{k}}}E_{I}\sum_{{\bf i}_{I}\in I_{2^{k}}}
  h^{2}I_{A_{k,l}}(X_{{\bf i}})>\gamma_{2^{k}}^{2}$. Hence
\[P(\sum_{{\bf i}\in I_{2^{k}}}h^{2}I_{A_{k,l}}
   (X_{{\bf i}})\geq \frac{1}{2}\gamma_{2^{k}}^{2})\geq 2^{-l-2}
  P(\exists_{i_{I^{'}}\in J_{2^{k}}}\ 2^{kl}E_{I}
  h^{2}I_{A_{k,l}}(\tilde{X}_{{\bf i}})> \gamma_{2^{k}}^{2})\]
and (\ref{ind2}) follows.  

(\ref{Bpr})$\Rightarrow$(\ref{Cpr}) This is the same as the above, except we use Lemma 
\ref{techdec} instead of Lemma \ref{techund}.
   
(\ref{C})$\Rightarrow$(\ref{A}). By the regularity assumptions (\ref{assd1}),
 (\ref{assd2}) and the Borel-Cantelli Lemma it is enough to prove that for
  any $t>0$
  \[\sum_{k=1}^{\infty}P(\frac{1}{\gamma_{2^{k}}}\max_{n\leq 2^{k}}
  |\sum_{{\bf i}\in I_{n}}\ve_{{\bf i}}h(X_{{\bf i}})|\geq t)<\infty.\]
  By our assumption (\ref{C}) it is enough to show that
  \[\sum_{k=1}^{\infty}P(\frac{1}{\gamma_{2^{k}}}\max_{n\leq 2^{k}}
    |\sum_{{\bf i}\in I_{n}}\ve_{{\bf i}}hI_{A_{k,d}}(X_{{\bf i}})|
    \geq t)<\infty.\]
  Since $d_{n}=\sum_{{\bf i}\in I_{n}}\ve_{{\bf i}}hI_{A_{k,d}}(X_{{\bf i}})$ 
  is a martingale, by Doob's maximal inequality we get
  \[P(\frac{1}{\gamma_{2^{k}}}\max_{n\leq 2^{k}}|\sum_{{\bf i}\in I_{n}}
    \ve_{{\bf i}}hI_{A_{k,d}}(X_{{\bf i}})|\geq t)\]
  \[\leq\frac{1}{t^{2}\gamma_{2^{k}}^{2}}E(\sum_{{\bf i}\in I_{2^{k}}}
  \ve_{{\bf i}}hI_{A_{k,d}}(X_{{\bf i}}))^{2}\leq
  \frac{2^{dk}}{t^{2}\gamma_{2^{k}}^{2}}Eh^{2}I_{A_{k,d}}(\tilde{X}).\]
  Thus it is enough to show that 
  \begin{equation}
  \label{aim1}
   \sum_{k=1}^{\infty}\frac{2^{dk}}{\gamma_{2^{k}}^{2}}
   Eh^{2}I_{A_{k,d}}(\tilde{X})<\infty.
  \end{equation}
  Let $\tau=\inf\{k: \tilde{X}\in A_{k,d}\}$, then 
  $\tilde{X}\in A_{\tau,d}\setminus A_{\tau-1,d}$ so by (\ref{assd3}) we get
  \[\sum_{k=1}^{\infty}\frac{2^{dk}}{\gamma_{2^{k}}^{2}}
   Eh^{2}I_{A_{k,d}}(\tilde{X})\leq E\sum_{k=\tau}^{\infty}
   \frac{2^{dk}}{\gamma_{2^{k}}^{2}}h^{2}(\tilde{X})\]
  \[\leq CE\frac{2^{d\tau}}{\gamma_{2^{\tau}}^{2}}h^{2}(\tilde{X})\leq
   C\sum_{k=1}^{\infty}E\frac{2^{dk}}{\gamma_{2^{k}}^{2}}h^{2}
   I_{A_{k,d}\setminus A_{k-1,d}}(\tilde{X}).\]
  Let us notice that by the definition of $A_{k,d}$ we have 
  $h^{2}(\tilde{X})I_{A_{k,d}\setminus A_{k-1,d}}(\tilde{X})\leq \gamma_{2^{k}}^{2}$ and
  $E_{I}2^{kl}h^{2}I_{A_{k,d}\setminus A_{k-1,d}}(\tilde{X})\leq 
  \gamma_{2^{k}}^{2}$
  for any $I\subset \{1,\ldots,d\}$ with $0<\card(I)=l<d$.
  Thus by Lemma \ref{techund}
  \[2^{-d-2}\min(\binom{2^{k-1}}{d}\frac{1}{\gamma_{2^{k}}^{2}}
    Eh^{2}I_{A_{k,d}\setminus A_{k-1,d}}(\tilde{X}),1)\leq 
    P(\sum_{{\bf i}\in I_{2^{k-1}}}h^{2}I_{A_{k,d}\setminus A_{k-1,d}}
    (X_{{\bf i}})>0)\]
   \[\leq P(\exists_{{\bf i}\in I_{2^{k-1}}}X_{{\bf i}}\in A_{k,d}\setminus A_{k-1,d})
   \leq P(\exists_{{\bf i}\in I_{2^{k-1}}} X_{{\bf i}}\notin A_{k-1,d}).\]
   So condition (\ref{C}) implies  that
   \[\sum_{k=1}^{\infty}\min(\frac{2^{dk}}{\gamma_{2^{k}}^{2}}
    Eh^{2}I_{A_{k,d}\setminus A_{k-1,d}}(\tilde{X}),1)<\infty\]
   and (\ref{aim1}) easily follows.
   
(\ref{Cpr})$\Rightarrow$(\ref{A}) and (\ref{Cpr})$\Rightarrow$(\ref{Apr}) In
the same way as above we show that (\ref{Cpr}) implies (\ref{aim1}) and that
(\ref{aim1}) implies (\ref{Apr}). $\Box$   

\bigskip

   The next Theorem will show how to deal with the condition (\ref{C}).
   Suppose that the sets $A_{k}$ are given and let us define the sets 
   $C_{k,l}$ and $B_{k,I}$ for $I\subset\{1,\ldots,d\}$ with
   $\card(I)=l$ by induction over $d-l$:
\[C_{k,d}=A_{k}\]
\[B_{k,I}=\{x_{I}\in E^{l}:2^{k(d-l)}\mu_{d-l}(C_{k,l+1}^{x_{I}})\geq 1 
  \mbox{ for } \card(I)=l\}\]
\[C_{k,l}=\{x\in C_{k,l+1}: x_{I}\notin B_{k,I}\mbox{ for all $I$ with } 
  \card(I)=l\}.\]
  
\begin{thm}
\label{Max}
\begin{equation}
\label{max1}
\sum_{k=1}^{\infty}P(\exists_{{\bf i}\in I_{2^{k}}}\ X_{{\bf i}}\in A_{k})
  <\infty
\end{equation}  
if and only if the following condition are satisfied
\begin{equation}
\label{max2}
  \forall_{l=1,\ldots,d-1} \forall_{I\subset\{1,\ldots,d\},\card(I)=l}\
  \sum_{k=1}^{\infty}P(\exists_{{\bf j}\in I_{2^{k}}^{l}}\ X_{{\bf j}}\in 
  B_{k,I})<\infty
\end{equation}  
\begin{equation}
\label{max3}
 \sum_{k=1}^{\infty}2^{kd}\mu_{d}(C_{k,1})<\infty. 
\end{equation} 
\end{thm}

\noindent{\bf Proof.} Let us notice that (\ref{max3}) immediately implies that
\[\sum_{k=1}^{\infty}P(\exists_{{\bf i}\in I_{2^{k}}}\ X_{{\bf i}}\in C_{k,1})
  <\infty.\]
Since by the definition of  sets $C_{k,l}$:
\begin{eqnarray*}\lefteqn{\{\exists_{{\bf i}\in I_{2^{k}}}\ X_{{\bf
i}}\in A_{k}\}}\\
& &\subset
  \{\exists_{{\bf i}\in I_{2^{k}}}\ X_{{\bf i}}\in C_{k,1}\}\cup
  \bigcup_{l=1}^{d-1}\bigcup_{I\subset\{1,\ldots d\},\card(I)=l}
 \{\exists_{{\bf i}_{I}\in I_{2^{k}}^{l}}\ X_{{\bf i}_{I}}\in B_{k,I}\},\end{eqnarray*}  
hence (\ref{max2}) and (\ref{max3}) imply (\ref{max1}).

To prove the second implication let us first notice that by the definition of
$C_{k,l}$ we have
\begin{equation}
\label{sect}
  2^{k(d-m)}\mu_{d-m}(C_{k,l}^{x_{I}})<1 \mbox{\ for any $I$ with }
  \card(I)=m\geq l.
\end{equation}  
Hence by Corollary 1  
\[P(\exists_{{\bf i}\in I_{2^{k}}}\ X_{{\bf i}}\in A_{k})\geq
 P(\exists_{{\bf i}\in I_{2^{k}}}\ X_{{\bf i}}\in C_{k,1})
\geq c_{d}2^{kd}\mu_{d}(C_{k,1})\]
so (\ref{max1}) implies (\ref{max3}). 

By Corollary 1 and (\ref{sect}) we also get that for any 
$I\subset \{1,\ldots,m\}$ with $\card(I)=l=1,\ldots,d-1$ we have
for $J=I^{c}$ and any $x_{I}\in E^{l}$
\[P(\exists_{{\bf i}_{J}\in I_{2^{k}}^{d-l}}\ X_{{\bf i}_{J}}\in 
  C_{k,l+1}^{x_{I}})\geq c_{d-l}2^{k(d-l)}\mu_{d-l}(C_{k,l+1}^{x_{I}}).\]
Thus
\[P(\exists_{{\bf i}\in I_{2^{k}}}\ X_{{\bf i}}\in A_{k})\geq
  P(\exists_{{\bf i}\in I_{2^{k}}}\ X_{{\bf i}}\in C_{k,l+1})
  \geq c_{d-l}P(\exists_{{\bf i}_{I}\in I_{2^{k}}^{l}}\ X_{{\bf i}_{I}}\in 
  B_{k,I})\]
and (\ref{max1}) implies (\ref{max2}). $\Box$

\subsection{Two-dimensional Case.}   

In the two-dimensional case let us define for $k=1,2,\ldots$ 
\begin{equation}
 \label{fk}
   f_{k}(x)= 2^{k}E_{Y}(h^{2}(x,Y)\wedge\gamma^{2}_{2^{k}}).
 \end{equation}
 
\begin{thm}
In the case of $d=2$ each of the equivalent conditions (\ref{A})-(\ref{Cpr})
is equivalent to the following condition
\begin{subequations}
\begin{equation}
\label{sub1}
\sum_{k=1}^{\infty}2^{k}P(f_{k}(X)\geq \gamma_{2^{k}}^{2})<\infty
\end{equation}
and
\begin{equation}
\label{sub2}
\sum_{k=1}^{\infty}2^{2k}P(h^{2}(X,Y)\geq \gamma_{2^{k}}^{2},
  f_{k}(X)< \gamma_{2^{k}}^{2}, f_{k}(Y)< \gamma_{2^{k}}^{2})<\infty.
\end{equation}
\end{subequations}
\end{thm}

\noindent{\bf Proof.} Again, we concentrate on the necessity, since the 
 sufficiency can be proved as in Theorem \ref{Sum2Max}. To obtain 
 (\ref{sub1}) first reduce to the decoupled sum of squares as in 
 Theorem \ref{Sum2Max} (\ref{Bpr}). One, then, has 
 \[P(\sum_{i,j\le 2^k}h^2(X_i,Y_j)\wedge{\gamma_{2^k}^2} > 
   \frac{1}{2}\gamma_{2^k}^2)\ge E_Y\max_{j\le 2^k}
   P_X(\sum_{i\le 2^k}h^2(X_i,Y_j)\wedge{\gamma_{2^k}^2} > 
   \frac{1}{2}\gamma_{2^k}^2)\]
 Applying Lemma \ref{techdec} (the case d=1) to the probability appearing 
in the last expectation, we see that 
 \[P_X(\sum_{i\le 2^k}h^2(X_i,Y_j)\wedge{\gamma_{2^k}^2} > 
   \frac{1}{2}\gamma_{2^k}^2)\ge\frac{1}{8}
   I_{\{2^kE_X(h^2\wedge{\gamma_{2^k}^2)} 
    > \gamma_{2^k}^2\}}\]
Hence, 
\begin{eqnarray*}
   \lefteqn{E_Y\max_{j\le 2^k}
  P_X(\sum_{i\le 2^k}h^2(X_i,Y_j)\wedge{\gamma_{2^k}^2} > 
  \frac{1}{2}\gamma_{2^k}^2)}\\
  & &\ge \frac{1}{8}
  P_Y(\max_{i\le 2^k}2^kE_X(h^2\wedge{\gamma_{2^k}^2}) 
  > \gamma_{2^k}^2)\\
  & &\ge\frac{1}{16}\min(1,2^kP_Y(2^{k}E_X(h^2\wedge\gamma_{2^k}^2)
  >\gamma_{2^k}^2)),
\end{eqnarray*}
 which implies (\ref{sub1}). But, we also have 
 \begin{eqnarray*}\lefteqn{P(\sum_{i,j\le 2^k}h^2(X_i,Y_j)\wedge
{\gamma_{2^k}^2}\ge \gamma_{2^k}^2)}\\
& &\ge P(\max_{i,j\le 2^k}h^2(X_i,Y_j)\wedge{\gamma_{2^k}^2}
I_{f_k(X_i),f_k(Y_j)\le \gamma_{2^k}^2}
\ge \gamma_{2^k}^2).
\end{eqnarray*}

Now, using the Section Lemma (Corollary 1) we have that the last quantity 
\[\ge 2^{-4}\min(1,2^{2k}P(h^2\wedge\gamma_{2^k}^2\ge \gamma_{2^k}^2, 
f_k(X), f_k(Y)< \gamma_{2^k}^2))
\]
And this implies (\ref{sub2}).$\Box$
\bigskip

\section{Convergence of series}

 In this section we will present the multidimensional generalizations of 
 symmetric case of Kolmogorov three series theorem, which states that for 
 independent random
 variables $X_{i}$ the following conditions are equivalent
\[\sum_{i=1}^{\infty} \ve_{i}X_{i} \mbox{ is a.s.\ convergent},\]
\[\sum_{i=1}^{\infty} X_{i}^{2}<\infty \mbox{ a.s.}\]
and
\[\sum_{i=1}^{\infty} E(X_{i}^{2}\wedge 1)< \infty.\]

 Let us first consider the two-dimensional case and define
\[c_{i}(x_{i})=\sum_{j=1}^{\infty}E_{Y}(h^{2}_{i,j}(x_{i},Y_{j})^{2}\wedge 1),\]
\[d_{j}(y_{j})=\sum_{i=1}^{\infty}E_{X}(h^{2}_{i,j}(X_{i},y_{j})^{2}\wedge 1).\]

\begin{thm}
  Suppose that the functions $c_{i},d_{j}$ are defined as above. Then the 
  following conditions 
  are equivalent
 \begin{equation}
 \label{ser2A}
   \lim_{n\rightarrow \infty}\sum_{i,j=1}^{n}\ve_{i}^{(1)}\ve_{j}^{(2)}
   h_{i,j}(X_{i},Y_{j}) \mbox{\ is a.s.\ convergent},
 \end{equation}
 \begin{equation}   
 \label{ser2B} 
   \sum_{i,j=1}^{\infty}h_{i,j}^{2}(X_{i},Y_{j})<\infty \mbox{\ a.s.}
 \end{equation}
and   
\begin{subequations} 
\label{ser2C}
\begin{equation}
\label{ser2C1}
c_{i}(X_{i})<\infty \mbox{ a.s. for all $i$ and } d_{j}(Y_{j})<\infty
\mbox{ a.s. for all $j$,}
\end{equation}
\begin{equation}
\label{ser2C2}
\sum_{i=1}^{\infty}P(c_{i}(X_{i})>1)<\infty \mbox{ and }
\sum_{j=1}^{\infty}P(d_{j}(Y_{j})>1)<\infty,
\end{equation}
\begin{equation}
\label{ser2C3}
\sum_{i,j=1}^{\infty}E(h_{i,j}^{2}(X_{i},Y_{j})\wedge 1)
  I_{\{c_{i}(X_{i})\leq 1, d_{j}(Y_{j})\leq 1\}}<\infty.
\end{equation}
\end{subequations}
\end{thm} 

\noindent{\bf Proof.} (\ref{ser2A})$\Leftrightarrow$(\ref{ser2B}). Let us first
 notice that 
 (\ref{ser2A}) and (\ref{ser2B}) are equivalent, respectively to the following
 two conditions
 \begin{equation}
 \label{ser2Apr}
 \forall_{\ve>0}\exists_{n}\ P(\sup_{k\geq n}
 |\sum_{n\leq i\vee j\leq k}
 \ve_{i}^{(1)}\ve_{j}^{(2)}h_{i,j}(X_{i},Y_{j})|>\ve)<\ve
 \end{equation}
 and
 \begin{equation}
 \label{ser2Bpr}
 \forall_{\ve>0}\exists_{n}\ P(\sum_{n\leq i\vee j}
 h_{i,j}^{2}(X_{i},Y_{j})>\ve)<\ve.
 \end{equation}
 By the hypercontractivity of Walshes (i.e., for sums of products of 
 Rademacher r.v.'s \cite{B} or \cite{KW}, sect. 3.4.) and the Paley-Zygmund 
 inequality we have
 \[P((\sum_{n\leq i\vee j}\ve_{i}^{(1)}\ve_{j}^{(2)}h_{i,j}(X_{i},Y_{j}))^{2}
   \geq\frac{1}{2}\sum_{n\leq i\vee j}h_{i,j}^{2}(X_{i},Y_{j}))\geq 
   \frac{1}{324}.\]
 Hence (\ref{ser2Apr}) implies (\ref{ser2Bpr}). On the other hand since 
 $d_{k}=\sum_{n\leq i\vee j\leq k}\ve_{i}^{(1)}\ve_{j}^{(2)}h(X_{i},Y_{j})$ is
 a martingale,
 we get by Doob's inequality
 \[P(\sup_{k\geq n}|\sum_{n\leq i\vee j\leq k}
  \ve_{i}^{(1)}\ve_{j}^{(2)}h_{i,j}(X_{i},Y_{j})|\geq t(\sum_{n\leq i\vee j}
  h_{i,j}^{2}(X_{i},Y_{j}))^{1/2})\leq t^{-2}\]
 and (\ref{ser2Bpr}) implies (\ref{ser2Apr}).

(\ref{ser2C})$\Rightarrow$(\ref{ser2B}). By condition (\ref{ser2C1}) we get
  that $\sum_{j=1}^{\infty}h_{i,j}^{2}(X_{i},Y_{j})<\infty$ a.s. for any $i$ and
  $\sum_{i=1}^{\infty}h_{i,j}^{2}(X_{i},Y_{j})<\infty$ a.s. for any $j$. Hence
  by condition (\ref{ser2C2}) it is enough to prove that
 \begin{equation}
 \label{aimser}
   Z=\sum_{i,j=1}^{\infty}(h_{i,j}^{2}(X_{i},Y_{j})\wedge 1)
   I_{\{c_{i}(X_{i})\leq 1, d_{j}(Y_{j})\leq 1\}}<\infty \mbox{ a.s.}.
 \end{equation}  
  However by Chebyshev's inequality 
 \[P(Z\geq t)\leq t^{-2}\sum_{i,j=1}^{\infty}E(h_{i,j}^{2}(X_{i},Y_{j})\wedge 1)
   I_{\{c_{i}(X_{i})\leq 1, d_{j}(Y_{j})\leq 1\}}\]
  and (\ref{aimser}) follows by (\ref{ser2C3}).    

(\ref{ser2B})$\Rightarrow$(\ref{ser2C}). Condition $c_{i}(X_{i})<\infty$ a.s.
  is equivalent to $\sum_{j=1}^{\infty}h_{i,j}^{2}(X_{i},Y_{j})<\infty$ a.s.,
  thus (\ref{ser2C1}) immediately follows by (\ref{ser2B}).

  To prove the condition (\ref{ser2C2}) let us notice that for sufficiently
  large $n$ we have
 \[P(\sum_{i=n,j=1}^{\infty}h_{i,j}^{2}(X_{i},Y_{j})\geq 
   \frac{1}{2})\leq 2^{-4}.\]
  Let us notice that by Lemma 1 (case $d=1$) we have for any $k\geq n$
 \[P_{Y}(\sum_{i=n,j=1}^{\infty}h_{i,j}^{2}(X_{i},Y_{j})\geq c_{k}(X_{k}))
   \geq
   P_{Y}(\sum_{j=1}^{\infty}h_{k,j}^{2}(X_{k},Y_{j})\wedge 1\geq c_{k}(X_{k}))
   \geq 2^{-3}.\]
  Thus
 \[P(\sum_{i=n,j=1}^{\infty}h_{i,j}^{2}(X_{i},Y_{j})\geq 
   \frac{1}{2})\geq 2^{-3}P(\max_{i\geq n}c_{i}(X_{i})>1),\]
  so $P(\max_{i\geq n}c_{i}(X_{i})>1)\leq 1/2$, which implies that
  $\sum_{i=1}^{\infty}P(c_{i}(X_{i}>1)<\infty$. In an analogous way
  we prove that $\sum_{j=1}^{\infty}P(d_{j}(Y_{j})>1)<\infty$.  

  Finally let
 \[m=\sum_{i,j=1}^{\infty}E(h_{i,j}^{2}(X_{i},Y_{j})\wedge 1)
   I_{\{c_{i}(X_{i})\leq 1, d_{j}(Y_{j})\leq 1\}}.\]
  We have
 \begin{eqnarray*}\lefteqn{E_{X}\sum_{i=1}^{\infty}h_{i,j}^{2}(X_{i},Y_{j})
   \wedge 1
   I_{\{c_{i}(X_{i})\leq 1, d_{j}(Y_{j})\leq 1\}}}\\
   & &\leq
   (E_{X}\sum_{i=1}^{\infty}h_{i,j}^{2}(X_{i},Y_{j})\wedge 1)
   I_{\{d_{j}(Y_{j})\leq 1\}}\leq 1
 \end{eqnarray*}
  and by a similar argument
 \[E_{Y}\sum_{j=1}^{\infty}(h_{i,j}^{2}(X_{i},Y_{j})\wedge 1)
   I_{\{c_{i}(X_{i})\leq 1, d_{j}(Y_{j})\leq 1\}}\leq 1.\]
  Hence by Lemma 1 we get
 \[P(\sum_{i,j=1}^{\infty}(h_{i,j}^{2}(X_{i},Y_{j})\wedge 1)
   I_{\{c_{i}(X_{i})\leq 1, d_{j}(Y_{j})\}}\geq 
   \frac{1}{2}m)\geq 2^{-4}\min(m,1),\]
  which implies that $m<\infty$. $\Box$

\bigskip

 Before formulating the result in the d-dimensional case we will need a 
few more definitions. 
 Let us define in this case $A_{0,{\bf i}}=E^{d}$ and then inductively 
 for $l=1,\ldots,d-1$, 
 $I\subset \{1,2,\ldots,d\}$ with $\card(I)=l$
  \[c_{{\bf i}_{I}}(x_{{\bf i}_{I}})=\sum_{{\bf i}_{I^{'}}}
   E_{I}^{'}(h^{2}_{({\bf i}_{I},{\bf i}_{I^{'}})}
   I_{A_{l-1,({\bf i}_{I},{\bf i}_{I^{'}})}}
   (x_{{\bf i}_{I}},\tilde{X}_{{\bf i}_{I^{'}}})\wedge 1),\]
 \[A_{l,{\bf i}}=\{x_{{\bf i}}\in A_{l-1,{\bf i}}:
   c_{{\bf i}_{I}}(x_{{\bf i}_{I}})\leq 1
   \mbox{\ for all $I$ with $\card(I)=l$}\}\]

\begin{thm}
Suppose that $c_{i_{I}}$ and $A_{l,{\bf i}}$ are defined as above.Then the following conditions 
are equivalent
\begin{equation}
\label{serdA}
 \sum_{{\bf i}\in Z_{+}^{d}}\ve_{{\bf i}}h_{{\bf i}}(\tilde{X}_{{\bf i}})
   \mbox{\ is a.s.\ convergent,}
\end{equation}   
\begin{equation}
\label{serdB}
 \sum_{{\bf i}\in Z_{+}^{d}}h^{2}_{{\bf i}}(\tilde{X}_{{\bf i}})
  <\infty \mbox{\ a.s.}
\end{equation}
and  
\begin{subequations}
\begin{equation} 
\label{serdC1}
 \sum_{{\bf i}_{I}\in Z_{+}^{d-1}}h^{2}_{{\bf i}}(\tilde{X}_{{\bf i}})<\infty
  \mbox{\ a.s.\ for all $I$ with $\card(I)=d-1$}
\end{equation}
\begin{equation}
\label{serdC2}
\sum_{{\bf i}_{I}\in Z_{+}^{l}}I_{\{c_{i_{I}}(\tilde{X}_{i_{I}})>1\}}<\infty
  \mbox{\ a.s.\ for all $I$ with $l=\card(I)=1,2\ldots,d-1$}
\end{equation}
\begin{equation}
\label{serdC3}  
  \sum_{{\bf i}\in Z_{+}^{d}}E(h_{{\bf i}}^{2}(\tilde{X}_{{\bf i}})\wedge 1)
  I_{A_{d-1,{\bf i}}}(\tilde{X}_{{\bf i}})<\infty
\end{equation}  
\end{subequations}
\end{thm} 

\noindent{\bf Proof.} As above.$\Box$

{\sc 
\noindent R.\ Lata{\l}a\\
  Institute of Mathematics\\
  Warsaw University\\
  Banacha 2\\
  02-097 Warszawa\\
  Poland}\\
email: {\tt rlatala@mimuw.edu.pl}

\bigskip

{\sc 
\noindent J.\ Zinn\\
 Deparment of Mathematics and Statistics\\
 Texas A\&M University\\
 College Station, Texas 77843}\\
email: {\tt jzinn@math.tamu.edu}

\end{document}